\documentclass[11pt,a4paper]{amsart}
\usepackage[all]{xy}
\usepackage{amscd,amssymb}

\newtheorem{theorem}{Theorem}[section]
\newtheorem{corol}[theorem]{Corollary}
\newtheorem{prop}[theorem]{Proposition}

\theoremstyle{definition}
\newtheorem{defin}[theorem]{Definition}

\theoremstyle{remark}
\newtheorem*{erem}{Remark}

 \def\larr{\longrightarrow}
 \def\={\setminus}	\def\ti{\tilde}
 \def\ver{\mathop\mathrm{Ver}\nolimits}
 \def\ind{\mathop\mathrm{ind}\nolimits}
\def\md#1{#1\mbox{-}\mathrm{mod}}
\def\vec{\mathsf{vec}}	
\def\bp{\bullet}	\def\iso{\simeq}
\def\arr{\mathop\mathrm{Arr}\nolimits}
\def\dar{\dashrightarrow}
\def\add{\mathop\mathrm{add}\nolimits}
\def\Fun{\mathop\mathrm{Fun}\nolimits}
\def\Rep{\mathop\mathrm{Rep}\nolimits}
\def\rep{\mathop\mathrm{rep}\nolimits}
\def\supp{\mathop\mathrm{supp}\nolimits}

\def\gdim{\mathop\mathrm{gl.dim}}
\def\Dim{\mathop\mathbf{dim}}
\def\iff{if and only if }
\def\oc{one-to-one correspondence }

\def\set#1{\left\{\,#1\,\right\}}
\def\setsuch#1#2{\left\{\,#1\,|\,#2\,\right\}}
\def\lst#1#2{ #1_1 , #1_2 , \dots , #1_{#2} }

\def\rad{\mathop\mathrm{rad}}

\def\ob{\mathop\mathrm{ob}\nolimits}
\def\Hom{\mathop\mathrm{Hom}}

\def\Ker{\mathop\mathrm{Ker}}
\def\im{\mathop\mathrm{Im}}
\def\8{\infty}		\def\+{\oplus}
\def\*{\otimes}	\def\dd{\partial}	

\def\al{\alpha}	\def\be{\beta}	\def\ga{\gamma}
\def\Ga{\Gamma}			\def\la{\lambda}
\def\La{\Lambda}	\def\si{\sigma}	\def\om{\omega}	
\def\De{\Delta}	\def\Si{\Sigma}	\def\eps{\varepsilon}
\def\tmu{\tilde\mu}

\def\bA{\mathbf A}	
\def\fM{\mathbf m}  	
\def\fR{\mathbf r}	\def\fB{\mathbf b}	
\def\fD{\mathbf d}	

\def\mN{\mathbb N}	\def\mZ{\mathbb Z}	
	
\def\kA{\mathcal A}		\def\kB{\mathcal B}
\def\kC{\mathcal C}	\def\kD{\mathcal D}
\def\kK{\mathcal K}	\def\kV{\mathcal V}		
\def\kW{\mathcal W}	\def\kQ{\mathcal Q}	
\def\kP{\mathcal P}		\def\kS{\mathcal S}		
\def\kM{\mathcal M}	\def\kU{\mathfrak U}	
\def\kJ{\mathcal J}

\def\dA{\mathfrak A}	\def\dB{\mathfrak B}
	\def\dP{\mathfrak P}
\def\dC{\mathfrak C}	

\def\sV{\mathsf V}		
\def\sI{\mathsf I}		\def\sJ{\mathsf J}
\def\sA{\mathsf A}		\def\sB{\mathsf B}
\def\sL{\mathsf L}		\def\sT{\mathsf T}
\def\sR{\mathsf R}

\def\tW{\tilde{\mathcal W}}	\def\oV{\overline{\mathcal V}}
\def\tA{\tilde{\mathcal A}}		\def\tV{\tilde{\mathcal V}}
\def\tB{\tilde{\mathcal B}}		\def\tU{\tilde{\mathcal U}}

\def\Mk{\Bbbk}

\begin{document}
 \author{Viktor I. Bekkert \and Yuriy A. Drozd}
 \thanks{The first author was supported by FAPESP (Grant N 98/14538-0) and CNPq (Grant 301183/00-7). }
 \title{Tame--wild dichotomy for derived categories}
\address{Departamento de Matem\'atica, ICEx, Universidade Federal de Minas 
Gerais, Av.  Ant\^onio Carlos, 6627, CP 702, CEP 30123-970, Belo 
Horizonte-MG, Brasil}
 \email{bekkert@mat.ufmg.br}
 \address{Department of Mechanics and Mathematics, Kyiv Taras Shevchenko University, 01033 Kyiv, Ukraine}
 \email{yuriy@drozd.org}
 \urladdr{http://drozd.org/\~{}yuriy}
 \subjclass[2000]{Primary: 16G60, secondary: 15A21, 16D90, 16E05}
  \begin{abstract}
  We prove that every finite dimensional algebra over an algebraically closed field is either derived tame
 or derived wild. The proof is based on the technique of boxes and reduction algorithm. It implies,
 in particular, that any degeneration of a derived wild algebra is derived wild; respectively, any
 deformation of a derived tame algebra is derived tame.
 \end{abstract} 
 \maketitle

 \section*{Introduction}

 The notions of tame and wild problems is now rather poplar in various branches of representation theory and related topics,
 especially because of the so-called \emph{tame-wild dichotomy} (cf. e.g. \cite{d0,dg} and other papers). 
 Namely, in most cases it so happens that either indecomposable representations depend on at most one
 parameter or their description becomes in some sense ``universal,'' i.e. containing a classification of representations
 of all finitely generated algebras. Last time these notions have also been studied for derived categories, and
 tame-wild dichotomy has been proved in some cases (though rather restrictive ones), cf. \cite{gk,br,ge}. 
 In this paper we shall prove such a dichotomy for derived categories of arbitrary finite dimensional algebras
 over an algebraically closed field. The used technique, just as in \cite{d0,dg} (see also the survey \cite{d1}),
 is that of ``matrix problems,'' more precisely, boxes and reduction algorithm. There are some new features: 
 we have to consider such boxes that the underlying category is no more free. Fortunately, the arising relations are of
 rather special nature, which leads to the notion of \emph{sliced boxes}. Actually, the tame-wild dichotomy is
 proved for such boxes, wherefrom the result for derived categories is obtained almost in the same way as the
 tame-wild dichotomy for representations of algebras has been obtained from that for free boxes. As it is rather
 usual, we formulate the result for \emph{locally finite dimensional categories}. If such a category only has finitely
 many indecomposable objects, this language is equivalent to that of finite dimensional algebras, though a bit more convenient.
 But categories with infinitely many indecomposables naturally arise in representation theory (for instance, when we consider
 coverings), so we prefer to use this language, especially as this generality does not imply the proofs.

  \section{Derived categories}
 \label{s1}

 We consider categories and algebras over a fixed algebraically closed field $\Mk$. A $\Mk$-category $\kA$ is called
 \emph{locally finite dimensional} (shortly \emph{lofd}) if the following conditions hold:
 \begin{enumerate}
\item  All spaces $\kA(x,y)$ are finite dimensional for all objects $x,y$.
 \item
  $\kA$ is \emph{fully additive}, i.e. it is additive and all idempotents in it split.\\
  Conditions 1,2 imply that the category $\kA$ is \emph{Krull--Schmidt},  i.e. each object uniquely decomposes into
 a direct sum of indecomposable objects; moreover, it is \emph{local}, i.e. for each indecomposable object $x$ the algebra
 $\kA(x,x)$ is local. We denote by $\ind\kA$ a set of representatives of isomorphism classes of indecomposable objects
 from $\kA$.
 \item
  For each object $x$ the set $\setsuch{y\in\ind\kA}{\kA(x,y)\ne0\text{ or }\kA(y,x)\ne0}$ is finite.
\end{enumerate}
 We denote by $\vec$ the category of finite dimensional vector spaces over $\Mk$ and by $\md\kA$ the category of 
 \emph{finite dimensional $\kA$-modules}, i.e. functors $M:\kA\to\vec$ such that $\setsuch{x\in\ind\kA}{Mx\ne0}$
 is finite. We also denote by $D(\kA)$ the derived category of the category $\md\kA$ and by $D^b(\kA)$ its full subcategory
 consisting of bounded complexes. The latter is again a lofd category.

 For an arbitrary category $\kC$ we denote by $\add\kC$ the minimal fully additive category containing $\kC$. For instance,
 one can consider $\add\kC$ as the category of finitely generated projective $\kC$-modules; especially, $\add\Mk=\vec$.
 We denote by $\Rep(\kA,\kC)$ the category of functors $\Fun(\kA,\add\kC)$ and call them \emph{representations} of
 the category $\kA$ in a category $\kC$. Obviously, $\Rep(\kA,\kC)\iso\Rep(\add\kA,\kC)$. If the category $\kA$ is lofd,
 we denote by $\rep(\kA,\kC)$ the full subcategory of $\Rep(\kA,\kC)$ consisting of the representations $M$ with \emph{finite
 support} $\supp M=\setsuch{x\in\ind\kA}{Mx\ne0}$. In particular, $\rep(\kA,\Mk)=\md\kA$.

 We denote by $\kD(\bA)$ (respectively, $\kD^-(\bA),\,\kD^b(\bA)\,$) the \emph{derived category} (respectively,
 right bounded and (two-sided) bounded derived category) of the category $\kA$-mod, where $\kA$ is a lofd category.
 Recall that $\kA$ embeds as a full subcategory into $\md\kA$. Namely, each object $x$ corresponds to the functor
 $\kA^x=\kA(x,\_\,)$. These functors are projective in the category $\md\kA$; if $\kA$ is fully additive, these are
 all projectives (up to isomorphism). On the other hand,
 $\md\kA$ embeds as a full subcategory into $\kD^b(\kA)$: a module $M$ is treated as a complex only having a
 unique nonzero component equal $M$ at the $0$-th position. It is also known that $\kD^-(\kA)$ can be identified
 with the category $\kK^-(\kA)$ whose objects are right bounded complexes of projective modules and morphisms
 are homomorphisms of complexes modulo homotopy \cite{gm}. If $\gdim\kA<\8$, every bounded complex has
 a bounded projective resolution, hence $\kD^b(\kA)$ can identified with $\kK^b(\bA)$, the category of
 bounded projective  complexes modulo homotopy, but it is no more the case if $\gdim\kA=\8$.
 Moreover, if $\kA$ is lofd, we can confine the considered complexes by \emph{minimal} ones, i.e.
 always suppose that $\im d_n\subseteq\rad P_{n-1}$ for all $n$.
 We denote by $\sJ$ the \emph{radical} of the category $\kA$, i.e. the set of morphisms having no invertible
 components with respect to some (hence, any) decomposition of its source and target into direct
 sums of indecomposables. Then $\rad M=\sJ M$ for each $M\in\md\kA$. 

 Even if $\gdim\bA=\8$, one easily shows \cite{d2} that $\kD^b(\kA)$ can be
 identified with a direct limit $\varinjlim_N\kQ^N(\kA)$ of the categories $\kQ^N(\kA)$ defined as follows.
 \begin{enumerate}
 \item  Objects of $\kQ^N(\kA)$ are right bounded complexes $(P_\bp,d_\bp)$ of projective modules from $\md\kA$
 with $P_n=0$ for $n>N$.
 \item
  Morphisms of $\kQ^N(\kA)$ are homomorphisms of complexes modulo \emph{quasi-homotopy}, where two homomorphisms
 $f,g:(P_\bp,d_\bp)\to(P'_\bp,d'_\bp)$ are said to be \emph{quasi-homotopic} if there are homomorphisms of modules
 $s_n:P_n\to P'_{n+1}$ such that $f_n-g_n=d'_{n+1}s_n+s_{n-1}d_n$ for all $n<N$.
 \item
  The functor $\kQ^N(\kA)\to \kQ^{N+1}(\kA)$ maps a complex $(P_\bp,d_\bp)$ to a complex  
 $$ 
   0\to \hat P_{n+1}\stackrel h\larr P_n\to P_{n-1}\to \dots \to P_m\to 0,
 $$ 
 where $h$ maps $\hat P_{n+1}$ onto $\Ker d_n$.
 (Such a complex is defined up to an isomorphism inside $\kQ^{N+1}(\kA)$.)
\end{enumerate}
 Note that these functors are full embeddings; thus all functors $\kQ^N(\kA)\to\kD^b(\kA)$ are full embeddings too,
 so we may treat $\kD^b(\kA)$ as a sort of union $\bigcup_N\kQ^N(\kA)$.
 Especially, in all classification problems, we may replace the study of the category $\kD^b(\kA)$ by that of the
 categories $\kQ^N(\kA)$.
 If $\kA$ is lofd, any complex from $\kQ^N(\kA)$ is isomorphic (in this category) to a
 minimal complex $P_\bp$ such that $\Ker d_N\subseteq\rad P_N$.
 We denote by $\kQ^N_0(\kA)$ the full subcategory of $\kQ^N(\kA)$ only consisting of such complexes.
 Thus $\kD^b(\kA)\iso\varinjlim_N\kQ^N_0(\kA)$.

  \begin{prop}\label{QN}
 Two complexes from $\kQ^N_0(\kA)$ are isomorphic in $\kQ^N(\kA)$  \iff they are isomorphic as complexes.  
 \end{prop} 
 \begin{proof}
 If two complexes $P_\bp,\,P'_\bp$ from $\kQ_0^N(\kA)$ are isomorphic in $\kQ^N(\kA)$, there is a diagram 
 $$ 
   \xymatrix{
   {P_N} \ar[rr]^{d_N} \ar@<.5ex>[d]^{\phi_N} &&
		 {P_{N-1}} \ar[rr]^{d_{N-1}} \ar@<.5ex>[d]^{\phi_{N-1}} &&
		 {P_{N-2}} \ar[r] \ar@<.5ex>[d]^{\phi_{N-2}} & {\dots} \\
   {P'_N} \ar[rr]^{d'_N} \ar@<.5ex>[u]^{\psi_N} &&
		 {P'_{N-1}} \ar[rr]^{d'_{N-1}} \ar@<.5ex>[u]^{\psi_{N-1}} &&
		 {P'_{N-2}} \ar[r] \ar@<.5ex>[u]^{\psi_{N-2}} & {\dots} 
	\ ,}
 $$ 
 where all upgoing and downgoing squares commute. Moreover, all products $\psi_n\phi_n\ (n<N)$ are of the form
 $1+\si_{n-1} d_n+d_{n+1}\si_n$, thus isomorphisms, as well as all products $\phi_n\psi_n\ (n<N)$.
 Hence all $\phi_n,\,\psi_n\ (n<N)$ are isomorphisms. 
 As $\phi_{N-1}(\im d_N)\subseteq\im d'_N$ and $\psi_{N-1}(\im d'_N)\subseteq\im d_N$, it implies that
 $\im d_N\iso\im d'_N$.
 Since $\Ker d_N\subseteq\rad P_N$, the latter is a projective cover of $\im d_N$, and $P'_N$ is a 
 projective cover of $\im d'_N$.
 Therefore $P_N\iso P'_N$.
 Moreover, $\phi_Nd'_N=\phi_{N-1}d_N:P_N\to\im d'_{N-1}$ is an epimorphism, hence
 $\im\phi_N+\Ker d'_N=P'_N$, so $\phi_N$ is an epimorphism, thus an isomorphism.
 \end{proof} 

 We introduce the notions of derived tame and derived wild lofd categories in the following way, which do not
 formally coincide with those of some earlier papers, such as \cite{br,ge,gk}, but is equivalent to them.
 Due to the preceding considerations, it is more convenient to deal with.
\newpage
  \begin{defin}\label{tw}
 Let $\kA$ be a lofd category. 
  \begin{enumerate}
\item   The \emph{rank} of an object $x\in\kA$ (or of the corresponding projective module $\kA^x$) is the function
 $\fR(x):\ind\kA\to\mZ$ such that $x\iso\bigoplus_{y\in\ind\kA}\fR(x)(y)y$. The \emph{vector rank} $\fR_\bp(P_\bp)$
 of a bounded complex of projective $\kA$-modules is the sequence $(\dots,\fR(P_n),\fR(P_{n-1}),\dots)$ (actually it
 only has finitely many nonzero entries).
 \item
  We call a \emph{rational family} of bounded minimal complexes over $\kA$ a bounded complex $(P_\bp,d_\bp)$
 of finitely generated projective $\kA\*\sR$-modules, where $\sR$ is a \emph{rational algebra},
 i.e. $\sR=\Mk[t,f(t)^{-1}]$ for a nonzero polynomial $f(t)$, and $\im d_n\subseteq\sJ P_{n-1}$
 For such a complex   we define $P_\bp(m,\la)$, where $m\in\mN,\,\la\in\Mk,\,f(\la)\ne0$, the complex
 $(P_\bp\*_\sR\sR/(t-\la)^m,d_\bp\*1)$. It is indeed a complex of projective $\kA$-modules. We put $\fR_\bp(P_\bp)=
 \fR_\bp(P_\bp(1,\la))$ (this vector rank does not depend on $\la$).
 \item
  We call a lofd category $\kA$ \emph{derived tame} if there is a set $\dP$ of rational families of bounded complexes over
 $\kA$ such that:
  \begin{enumerate}
 \item  For each vector rank $\fR_\bp$ the set $\dP(\fR_\bp)=\setsuch{P_\bp\in\dP}{\fR_\bp(P_\bp)=\fR}$ is finite.
 \item
  For each vector rank $\fR_\bp$ all indecomposable complexes $(P_\bp,d_\bp)$ of projective $\kA$-modules 
 of this vector rank, except finitely many isomorphism classes, are isomorphic to $P_\bp(m,\la)$ for some $P_\bp\in\dP$
 and some $m,\la$.
\end{enumerate}
 The set $\dP$ is called a \emph{parameterising set} of $\kA$-complexes.
 \item
  We call a lofd category $\kA$ \emph{derived wild} if there is a bounded complex $P_\bp$ of projective modules over
$\kA\*\Si$, where $\Si$ is the free $\Mk$-algebra in 2 variables, such that, for every finite dimensional $\Si$-modules
 $L,L'$,
  \begin{enumerate}
\item   $P_\bp\*_\Si L\iso P_\bp\*_\Si L'$ \iff $L\iso L'$.
 \item
  $P_\bp\*_\Si L$ is indecomposable \iff so is $L$.
\end{enumerate}
 (It is well-known that then an analogous complex of $\kA\*\Ga$-modules exists for every finitely generated $\Mk$-algebra
 $\Ga$.)
\end{enumerate} 
 \end{defin} 
 
 Note that, according to these definitions, every \emph{derived discrete} (in particular, \emph{derived finite}) lofd category
 \cite{vo} is derived tame (with the empty set  $\dP$).
 Simple geometric considerations, like in \cite{d3}, show that neither lofd category can be both derived tame and 
 derived wild. We are going to demonstrate the following result.

   \begin{theorem}[\sc Main Theorem]\label{main}
  Every lofd category over an algebraically closed field is either derived tame or derived wild. 
 \end{theorem} 

 This theorem will be proved in Section \ref{s3}. Note that, in particular, it makes valid the following corollaries,
 which have been proved in \cite{d2} under supposition that every finite dimensional algebra is either derived tame or derived wild.

   \begin{corol}\label{12}
  Let $\kA$ be a flat family of finite dimensional algebras based on an algebraic variety $X$. Then
 the set $\setsuch{x\in X}{\kA(x) \text{ \em is derived wild}}$ is a union of a countable sequence of closed subsets. 
 \end{corol} 

  \begin{corol}\label{13}
  Suppose that a finite dimensional algebra $\sA$ \emph{degenerates} to another algebra $\sB$ (or, the same,
 $\sB$ \emph{deforms} to $\sA$), i.e. there is a flat family of algebras $\kA$ based on a variety $X$ such
 that $\kA(x)\simeq\sA$ for all $x$ from a dense open subset $U\subseteq X$ and there is a point $y\in X$
 such that $\kA(y)\simeq\sB$. If $\sA$ is derived wild, so is $\sB$; respectively, if $\sB$ is derived tame,
 so is $\sA$. 
 \end{corol} 
 
 If the families are not assumed flat, these assertions are no more true \cite{br}
  (see \cite{d2,d4} for further comments).

 \section{Related boxes}
 \label{s2}

 Recall \cite{d0,d1} that a \emph{box} is a pair $\dA=(\kA,\kV)$ consisting of a category $\kA$ and an 
 $\kA$-coalgebra $\kV$. We denote by $\mu$ the comultiplication in $\kV$, by $\eps$ its counit and
 by $\oV=\Ker\eps$ its \emph{kernel}. We always suppose that $\dA$ is \emph{normal}, i.e. there is
 a \emph{section} $\om:x\mapsto\om_x\ (x\in\ob\kA)$ such that $\eps(\om_x)=1_x$ and $\mu(\om_x) 
 =\om_x\*\om_x$ for all $x$. A category $\kA$ is called \emph{free} if it is isomorphic to a path category
 $\Mk\Ga$ of an oriented graph (quiver) $\Ga$, and \emph{semi-free} if $\kA=\Mk\Ga[\sL^{-1}]$, where
 $\sL$ is a set of \emph{loops}, i.e. arrows $a:x\to x$ from $\Ga$. The arrows $a:x\to y$ with $x\ne y$
 will be called \emph{edges}. If $\Ga$ contains no arrows at all, the category $\kA$ is called \emph{trivial};
 if $\Ga$ only has loops, and at most one loop at every vertex $x$, $\kA$ is called \emph{minimal}. A normal
 box $\dA=(\kA,\kV)$ is called \emph{free} (\emph{semi-free}) if so is the category $\kA$, while the kernel
 $\oV$ is a free $\kA$-bimodule. If we fix a set of free generators $\De$ of $\oV$, we call the elements from $\De$
 \emph{dashed arrows} of the box $\dA$, while the arrows of $\Ga$ are called \emph{solid arrows}. The 
 union $\arr\dA=\Ga\cup\De $ is called a \emph{set of free} (or \emph{semi-free}) \emph{generators} of
 the free (semi-free) box $\dA$. We also call the objects of $\kA$ the \emph{vertices} of $\dA$, denote
 by $\ver\dA$ the set of vertices, and write $\arr^0\dA=\Ga,\ \arr^1\dA=\De$.
 Note that a choice of free (semi-free) generators is usually not unique, and most
 of proofs related to boxes use a change of free (semi-free) generators. 

 Recall that the \emph{differential} of a normal box $\dA=(\kA,\kV)$ is the pair $\dd=(\dd_0,\dd_1)$ of
 mappings, $\dd_0:\kA\to\oV,\ \dd_1:\oV\to\oV\*_\kA\oV$, namely
 \begin{align*}
  \dd_0 a&=a\om_x-\om_ya \quad\text{ for } a\in \kA(x,y),\\
  \dd_1 v&=\mu(v)-v\*\om_x-\om_y\*v \quad\text{for } v\in\oV(x,y).
 \end{align*} 
 Usually we omit the index and write both $\dd a$ and $\dd v$.
 A set of arrows $\arr\dA$ of semi-free box is said to be \emph{triangular}, if there is a function
 $h:\arr\dA\to\mN$ (called \emph{height}) such that, for any $a\in\arr\dA$ (either solid or dashed)
 $\dd a$ belongs to the sub-box generated by the arrows $b\in\arr\dA$ with $\dd b<\dd a$,
 especially, $\dd a=0$ if $h(a)=0$. If such a set of arrows exists, we call the box $\dA$ \emph{triangular}.

 A normal box $\dA$ such that the category $\kA$ is trivial, is called \emph{so-trivial} (trivial with respect to
 solid arrows). If $\kA$ is a minimal category and $\dd a=0$ for each solid loop, we call the box $\dA$
 \emph{so-minimal}.

 In what follows we also use boxes, which are not free (or semi-free), but are their factors. Namely, let $\dA=(\kA,\kV)$
 be a semi-free box, $\sI\subseteq\kA$ be an ideal of the category $\kA$ such that $\dd a\in\sI\oV+\oV\sI$ for all $a\in\sI$.
 Denote by $\dA/\sI$ the box $(\tA,\tV)$, where $\tA=\kA/\sI$ and $\tV=\kV/(\sI\kV+\kV\sI)$, with
 natural comultiplication and counit. Note that in this case the kernel of the box $\dA/\sI$ is a free $\tA$-bimodule,
 namely, it is isomorphic to $\oV/(\sI\oV+\oV\sI$). If $\dA$ is a triangular semi-free box and the ideal $\sI$ is contained in the ideal
 generated by all products $ab$, where $a,b$ are solid arrows, we call $\ti\dA=\dA/\sI$ a \emph{convenient} box. The
 vertices and arrows of $\ti\dA$ are, by definitions, those of $\dA$. Especially, the notions of \emph{triangular set of arrows}
 and \emph{triangular box} are transmitted to convenient boxes. Actually, we need a rather specific kind of
 convenient boxes, defined as follows.

  \begin{defin}\label{slice}
 \begin{enumerate}
\item  A free box $\dA$ is called \emph{sliced} if there is a function $s:\ver\dA\to\mZ$ such that
 \begin{enumerate}
 \item  $s(y)<s(x)$  for every solid arrow $a:x\to y$; we set $s(a)=s(x)$;
 \item
  $s(x)=s(y)$ for every dashed arrow $\ga:x\dar y$; we set $s(\ga)=s(x)=s(y)$.
\end{enumerate}
 \item
  A box $\ti\dA=\dA/\sI$, where $\dA=(\kA,\kV)$ is a free box and $\sI\subset\kA$ is an ideal in $\kA$
 such that $\dd a\in\sI\oV+\oV\sI$ for all $a\in\sI$, is called \emph{sliced}
 if so is the box $\dA$.
 \end{enumerate}
 We call the function $s$ a \emph{slicing} of the box $\dA$ or $\ti\dA$.
 \end{defin} 

 Note that if a free box $\dA$ is sliced, there are neither loops nor oriented cycles in it. Therefore, if an ideal $\sI$
 is not contained in the ideal generated by the paths of length $2$, we are able just drop an arrow that occur in an
 element of $\sI$. Hence sliced boxes are always convenient. 

 A \emph{representation} of a box $\dA=(\kA,\kV)$ over a category $\kC$ is defined as a functor
 $M:\kA\to\add\kC$. A \emph{morphism} of such representations $f:M\to N$ is defined as a 
 homomorphisms of $\kA$-modules $\kV\*_\kA M\to N$. If $g:N\to L$ is another morphism, there
 product is defined as the composition 
 $$ 
 \begin{CD} 
   \kV\*_\kA M @>\mu\*1>> \kV\*_\kA\kV\*_\kA M @>1\*f>>\kV\*_\kA N@>g>> L.
 \end{CD} 
 $$ 
 Thus we obtain the \emph{category of representations} $\Rep(\dA,\kC)$.
 If $\dA$  is a free (or a convenient) box, we denote by $\rep(\dA,\kC)$ the full subcategory of $\Rep(\dA,\kC)$
 consisting of representations with finite support $\supp M=\setsuch{x\in\ver\dA}{Mx\ne0}$.
 If $\kC=\vec$, we write $\Rep(\dA)$ and $\rep(\dA)$. 

 Given a lofd $\kA$, we are going to construct a sliced box $\dB=\dB(\kA)=(\kB,\kW)$ such that its 
 representations classify the objects of the derived category $\kD^b(\kA)$. 
 
 We denote by $\kS$ the trivial category with the set of objects 
 $$ 
   \ob\kS=\setsuch{(x,n)}{x\in\ind\kA,\,n\in\mZ}
 $$ 
  and consider the $\kS$-bimodule $\kJ$ such that  
 $$ 
   \kJ\big((x,n),(y,m)\big)= \begin{cases}
  0 &\text{if } m\ne n-1,\\
  \sJ(x,y)^* &\text{if }m=n-1,
 \end{cases} 
 $$ 
 where $\sJ$ is the radical of $\kA$ and $V^*$ denotes the dual vector space to $V$.
 Let $\tB=\kS[\kJ]$ be the tensor category of this bimodule; equivalently, it is  the free category having the same
 set of objects as $\kS$ and the union of bases of all $\kJ\big((x,n),(y,m)\big)$ as a set of free generators. 
 Denote by $\kU$ the $\kS$-bimodule such that 
  $$
 \kU\big((x,n),(y,m)\big)= \begin{cases}
  0 &\text{if } n\ne m,\\
 \kA(x,y)^* &\text{if } n=m
 \end{cases} 
 $$
 and set $\tU=\tB\*_\kS\kU\*_\kS\tB$. Dualizing the multiplication $\kA(y,z)\*\kA(x,y)\to\kA(x,z)$, we
 get homomorphisms
 \begin{align*}
  \la_r&: \tB\larr \tB\*_\kS\tU,\\
  \la_l&: \tB\larr \tU\*_\kS\tB,\\
  \tmu&: \tU\larr \tU\*_\kS\tU.
 \end{align*} 
  In particular, $\tmu$ defines on $\tU$ a structure of $\tB$-coalgebra. Moreover, the sub-bimodule $\kU_0$
 generated by $\im(\la_r-\la_l)$ is a coideal in $\tU$, i.e.
 $\tmu(\kU_0)\subseteq\kU_0\*_{\tB}\tU\+\tU\*_{\tB}\kU_0$.
 Therefore, $\tW=\tU/\kU_0$ is also a $\tB$-coalgebra, so we get a box $\ti\dB=(\tB,\tW)$. One easily checks,
 like in \cite{d0}, that it is free and triangular. 

 Dualizing multiplication also gives a mapping
 \begin{equation}\label{e21} 
 \nu:\sJ(x,y)^*\larr\bigoplus_z\sJ(z,y)^*\*\sJ(x,z)^*. 
 \end{equation} 
 Namely, if we choose bases $\set\al,\,\set\be\,\set\ga$ in the spaces, respectively, $\sJ(x,y),$ $\sJ(z,y),\,\sJ(x,z)$,
 and dual bases $\set{\al^*},\,\set{\be^*},\,\set{\ga^*}$ in their duals, then $\be^*\*\ga^*$ occurs in $\nu(\al^*)$
 with the same coefficient as $\al$ occurs in $\be\ga$.
 Note that the right-hand space in \eqref{e21} coincide with each $\tB\big((x,n),(y,n-2)\big)$. Let $\sI$ be the ideal
 in $\tB$ generated by the images of $\nu$ in all these spaces and $\dB=\ti\dB/\sI=(\kB,\kW)$, where
 $\kB=\tB/\sI,\ \kW=\tW/(\sI\tW+\tW\sI)$. One easily checks that $\dd\sI\subseteq\sI\tW+\tW\sI$, so it is a convenient box.
  If necessary, we write $\dB(\kA)$ to emphasise that this box has been 
 constructed from a given algebra $\kA$. Certainly, $\dB$ is a sliced triangular box, and the following result holds.

  \begin{theorem}\label{box}
  The category of finite dimensional representations $\rep(\dB(\kA))$ is equivalent to the category $\kP^b_{\min}(\kA)$
 of bounded minimal projective $\kA$-complexes.
 \end{theorem} 
 \begin{proof}
 We denote $\sJ^x=\sJ(x,\_\,)=\rad\kA^x$. Then $\Hom_\kA(\kA^x,\sJ^y)\simeq\sJ(x,y)$.
 A representation $M\in\rep(\dB)$ is given by vector spaces $M(x,n)$ and linear mappings 
 $$ 
 M_{xy}(n):\sJ(x,y)^*=\kA\big((x,n),(y,n-1)\big)\to\Hom\big(M(x,n),M(y,n-1)\big),
 $$
 where $x,y\in\ind\kA,\,n\in\mZ$, subject to the relations 
 \begin{equation}\label{e22}
   \sum_z \fM\big(M_{zy}(n)\*M_{xz}(n+1)\big)\nu(\al)=0
 \end{equation}
 for all $x,y,n$ and all $\al\in\sJ_{xy}$, where $\fM$ denotes the multiplication of mappings 
 \begin{multline*} 
   \Hom\big(M(z,n),M(y,n-1)\big)\*\Hom\big(M(x,n+1),M(z,n)\big)\to\\
	\to\Hom\big(M(x,n+1),M(y,n-1)\big).   
 \end{multline*}  
 For such a representation, set $P_n=\bigoplus_x \kA^x\*M(x,n)$. Then
 $\rad P_n=\bigoplus_x \sJ^x\*M(x,n)$ and
 \begin{align*}
  \Hom_\kA(P_n,\rad P_{n-1})&\simeq \bigoplus_{x,y} \Hom_\kA\big(\kA^x\*M(x,n),\sJ^y\*M(y,n-1)\big)\simeq\\
	&\simeq \bigoplus_{x,y} \Hom\big(M(x,n),\Hom_\kA\big(\kA^x,\sJ^y\*M(y,n-1)\big)\big)\simeq\\
	&\simeq \bigoplus_{x,y} M(x,n)^*\*\sJ(x,y)\*M(y,n-1) \simeq\\
	&\simeq \bigoplus_{x,y} \Hom\big(\sJ^*(x,y),\Hom\big(M(x,n),M(y,n-1)\big)\big).
 \end{align*} 
 Thus the set $\setsuch{M_{xy}(n)}{x,y\in\ind\kA}$ defines a homomorphism $d_n:P_n\to P_{n-1}$ and vice versa.
 Moreover, one easily verifies that the condition \eqref{e22} is equivalent to the relation $d_nd_{n+1}=0$.
 Since every projective $\kA$-module can be given in the form $\bigoplus_x\kA^x\*V_x$ for some
 uniquely defined vector spaces $V_x$, we get a \oc between finite dimensional representations of $\dB$
 and bounded minimal complexes of projective $\kA$-modules.
 In the same way one also establishes \oc between morphisms of representations and of the corresponding
 complexes, compatible with their multiplication, which accomplishes the proof.
  \end{proof} 

 Note that we can pick up subcategories of $\rep(\dB)$ that describe each of $\kQ^N(\kA)$.
 Namely, denote by $\rep^N(\dB)$ the full subcategory of $\rep(\dB)$ consisting of all representations
 $M$ such that $M(x,n)=0$ for $n>N$.
 Let $\sT_N$ be the ideal of $\rep^N(\dB)$ generated by the identity morphisms of
 \emph{trivial representations} $S_{x,N}$, where $S_{x,N}(x,N)=\Mk$, $S_{x,N}(y,n)=0$ if
 $(y,n)\ne(x,N)$.
 Obviously, the equivalence of the categories $\rep(\dB)$ and $\kP_{\min}^b(\kA)$ maps representations
 from $\rep^N(\dB)$ onto the complexes $P_\bp$ with $P_n=0$ for $n>N$. 
 Moreover, it maps $S_{x,N}$ to the complex $T^{x,N}_\bp$ with $T^{x,N}_N=\kA^x$, 
 $T^{x,N}_n=0$ for $n\ne N$.
 Note that a morphism of complexes from $\kQ^N(\kA)$ is quasi-homotopic to zero \iff it factorises through
 a direct sum of complexes $T^{x,N}_\bp$.
 It gives the following

 \begin{corol}\label{boxN}
  The category $\kQ^N(\kA)$ is equivalent to the factor category $\rep^N(\dB)/\sT_N$.
 \end{corol} 

\noindent
 Evidently, $\,\ind\big(\rep^N(\dB)/\sT_N\big)=
 \ind\big(\rep^N(\dB)\big)\=\setsuch{S_{x,N}}{x\in\ver\kA}$.

  \begin{corol}\label{twbox}
  An algebra $\kA$ is derived tame (derived wild) if so is the box $\dB(\kA)$. 
 \end{corol}

 \section{Proof of the Main Theorem}
 \label{s3}

 Now we are able to prove the main theorem. Namely, according to Corollary \ref{twbox}, it follows from the
 analogous result for sliced boxes.

  \begin{theorem}\label{mbox}
 Every sliced triangular box is either tame or wild.   
 \end{theorem} 

 Actually, just as in \cite{d0} (see also \cite{d1}), we shall prove this theorem in the following form.

\theoremstyle{plain}
\newtheorem*{31a}{Theorem 3.1a}
  \begin{31a}\label{mini}
  Suppose that a sliced triangular box $\dA=(\kA,\kV)$ is not wild. For every dimension $\fD$ of its
 representations there is a functor $F_\fD:\kA\to\add\kM$, where $\kM$ is a minimal category, such that
 every representation $M:\kA\to\vec$ of $\dA$ of dimension $\Dim (M)\le\fD$ is isomorphic to the
 inverse image $F^*N=N\circ F$ for some  functor $N:\kM\to\vec$. Moreover, $F$ can be chosen
 \emph{strict}, which means that $F^*N\simeq F^*N'$ implies $N\simeq N'$ and $F^*N$ is
 indecomposable if so is $N$.
 \end{31a} 

  \begin{erem}\label{r21}
  We can consider the induced box $\dA^F=(\kM,\kM\*_\kA\kV\*_\kA\kM)$. It is a so-minimal box,
 and $F^*$ defines a full and faithful functor $\rep(\dA^F)\to\rep(\dA)$. Its image consists of all
 representations  $M:\kA\to\vec$ that factorise through $F$. 
 \end{erem} 

 \begin{proof}
 As we only consider finite dimensional representations, we may assume that the set of
 objects is finite.
 Then we may assume that all values of a slicing $s:\ver\dA\to\mZ$ belong to $\mN$, and there are
 finitely many of them. 
 Let  $m=\max\setsuch{s(x)}{x\in\ver\dA}$.
 We use induction on $m$. 
 If $m=1$, $\dA$ is free, and our claim has been proved in \cite{d0}.
 So we may suppose that the theorem is true for smaller
 values of $m$, especially, it is true for the restriction $\dA'=(\kA',\kV')$ of the box $\dA$
 onto the subset $\sV=\setsuch{x\in\ver\dA}{s(x)<m}$.
 Thus there is a strict functor $F':\kA'\to\add\kM$,
 where $\kM$ is a minimal category, such that every representation of $\dA'$ of dimension smaller
 than $\fD$ is of the form ${F'}^*N$ for $N:\kM\to\vec$.
 Consider now the amalgamation $\kB=\kA\bigsqcup^{\kA'}\kM$ and the box
 $\dB=(\kB,\kW)$, where $\kW=\kB\*_\kA\kV\*_\kA\kB$.
 The functor $F'$ extends to a functor $F:\kA\to\kB$ and
 induces a homomorphism of $\kA$-bimodules $\kV\to\kW$; so it defines a functor $F^*:
 \rep(\dB)\to\rep(\dA)$, which is full and faithful.
 Moreover, every representation of $\dA$ of dimension smaller than $\fD$
 is isomorphic to $F^*N$ for some $N$, and all possible
 dimensions of such $N$ are restricted by some vector $\fB$.
 Therefore, it is enough to prove the claim of the theorem for the box $\dB$. 

 Note that the category $\kB$ is generated by the loops from $\kM$ and the images of arrows from
 $\kA(a,b)$ with $s(a)=m$ (we call them \emph{new arrows}).
 It implies that all possible relations between these morphisms are of the form
 $\sum_\al g_\al(\be)\al=0$, where $\be\in\kB(b,b)$ is a loop (necessarily minimal,
 i.e. with $\dd\be=0$), $g_\al$ are some polynomials, and $\al$ runs through the set
 of new arrows from $a$ to $b$ for some $a$ with $s(a)=m$.
 Consider all of these relations for a fixed $a$; let them be $\sum_\al g_{\al,k}(\be)\al=0$.
 Their coefficients form a matrix $\big(g_{\al,k}(\be)\big)$. Using transformations
 of the set $\set\be$ and of the set of relations, we can make this matrix diagonal,
 i.e. make all relations being $ f_\al(\be)\al=0$ for some polynomials $f_\al$.
 If one of $f_\al$ is zero, the box $\dB$ has a sub-box  
$$
 \xymatrix{
	{a}  \ar[rr]^{\al} && b \ar@(ur,dr)[]^{\be}	},
$$

 \medskip\noindent
 with $\dd\al=\dd\be=0$,
 which is wild; hence $\dB$ and $\dA$ are also wild. Otherwise, let $f(\be)\ne0$ be a common multiple
 of all $f_\al(\be)$, $\La=\set{\lst \la r}$ be the set of roots of $f(\be)$. If $N\in\rep(\dB)$ is such
 that $N(\be)$ has no eigenvalues from $\La$, then $f(N(\be))$ is invertible; thus $N(\al)=0$ for all
 $\al:a\to b$. So we can apply the \emph{reduction of the loop} $\be$ with respect to the set
 $\La$  and the dimension $d(a)$, as in \cite[Propositions 3,4]{d0} or \cite[Theorem 6.4]{d1}.
 It gives a new box that has the same number of loops as $\dB$,
 but the loop corresponding to $\be$ is ``isolated,'' i.e. there are no more arrows starting or ending
 at the same vertex. In the same way we are able to isolate all loops, obtaining a semi-free triangular box 
 $\dC$ and a morphism $G:\dB\to\dC$ such that $G^*$ is full and faithful and all representations of
 $\dB$ of dimensions smaller than $\fB$ are of the form $G^*L$. As the theorem is true for semi-free
 boxes, it accomplishes the proof.
 \end{proof}

 \end{document}